# RANDOM MATRIX CENTRAL LIMIT THEOREMS FOR NONINTERSECTING RANDOM WALKS

By Jinho Baik[1] and Toufic M. Suidan[2]

*University of Michigan, Ann Arbor and University of California, Santa Cruz*

We consider nonintersecting random walks satisfying the condition that the increments have a finite moment generating function. We prove that in a certain limiting regime where the number of walks and the number of time steps grow to infinity, several limiting distributions of the walks at the mid-time behave as the eigenvalues of random Hermitian matrices as the dimension of the matrices grows to infinity.

**1. Introduction.** It is known that various limiting local statistics arising in random matrix theory are independent of the precise structure of the randomness of the ensemble [11, 17, 18, 21, 32, 44, 48]. For example, consider the set of Hermitian matrices equipped with a probability measure invariant under unitary conjugation. For a very general class of measures, as the size of the matrix becomes large, the largest eigenvalue converges in distribution to the Tracy–Widom distribution, while the gap probability in the "bulk scaling limit" converges to a (different) universal distribution.

It has been discovered that the limiting distributions arising in random matrix theory also describe limit laws of a number of specific models in combinatorics, probability theory and statistical physics; apparently, these models are not expressible in terms of random matrix ensembles. Examples include the longest increasing subsequence of random permutations [5, 13, 31, 43], random Aztec and Hexagon tiling models [9, 33], last passage percolation models with geometric and exponential random

Received April 2006.
[1]Supported in part by NSF Grant DMS-04-57335 and the Sloan Fellowship.
[2]Supported in part by NSF Grants DMS-02-02530 and DMS-05-53403.

variables [30], polynuclear growth models [34, 45] and vicious walker models [3, 27]. For these models, the distribution function of interest was computed explicitly in terms of certain determinantal formulae and the asymptotic analysis of these determinants yielded the desired limit law. Nevertheless, it is believed that such limit laws should hold for a class of models much wider than the explicitly computable ("integrable") models. One such universality result for models "outside random matrices" was obtained in [10, 12, 50] for thin last passage percolation models with general random variables.

This paper studies nonintersecting random walks and proves random matrix central limit theorems in a certain limiting regime. The motivation for this study comes from two sources. The first is the fact that the eigenvalue density function of the Gaussian unitary ensemble can be described in terms of a nonintersecting Brownian bridge process [22, 33]. Namely, consider $n$ standard Brownian bridge processes $(B_t^{(1)},\dots,B_t^{(n)})$ conditioned not to intersect during the time interval $(0,2)$ (i.e., $B_t^{(1)} > \cdots > B_t^{(n)}$ for $0<t<2$), all starting from and ending at the origin. A simple computation shows that the distribution of $\{B_1^{(1)},\dots,B_1^{(n)}\}$ at time 1 is the same as the distribution of the eigenvalues of the $n\times n$ Gaussian unitary ensemble; see Section 1.1.1 below for the computation. Hence, it is natural to ask if the same limit laws hold for general nonintersecting random walks. The second motivation is that a number of the aforementioned specific probability models for which the random matrix central limit theorem was obtained are indeed interpreted in terms of nonintersecting random walks. We mention a few of them in the following subsection.

1.1. *Motivating examples.* We begin by introducing two distribution functions. Define the kernels

$$\mathbb{A}(a,b)=\frac{\mathrm{Ai}(a)\mathrm{Ai}'(b)-\mathrm{Ai}'(a)\mathrm{Ai}(b)}{a-b},\qquad \mathbb{S}(a,b)=\frac{\sin(\pi(a-b))}{\pi(a-b)}. \tag{1}$$

Set

$$F_{\mathrm{TW}}(\xi)=\det(1-\mathbb{A}|_{(\xi,\infty)}),\qquad F_{\mathrm{Sine}}(\eta)=\det(1-\mathbb{S}|_{[-\eta,\eta]}). \tag{2}$$

The Tracy–Widom distribution, $F_{\mathrm{TW}}$, is the limiting distribution of the largest eigenvalue and $F_{\mathrm{Sine}}$ is the limiting distribution for the gap probability of the eigenvalues "in bulk" in Hermitian random matrix theory.

1.1.1. *Nonintersecting Brownian bridge process.* Let $B_t=(B_t^{(1)},\dots,B_t^{(n)})$ be an $n$-dimensional standard Brownian motion. We compute the density function of $B_1$ conditioned that $B_t^{(1)}>B_t^{(2)}>\cdots>B_t^{(n)}$ for $0<t<2$ and $B_0 = B_2 = (0, \dots, 0)$. Let $p_t(x,y) = \frac{1}{\sqrt{2\pi t}} e^{-(x-y)^2/(2t)}$. The argument of Karlin and McGregor [36] implies that the density function of $n$ one-dimensional nonintersecting Brownian motions at time $t$ which start from $(x_1, \dots, x_n)$, where $x_1 > \cdots > x_n$, is given by

$$f_t(b_1, \dots, b_n) = \det(p_t(x_i, b_j))_{i,j=1}^n, \qquad b_1 > \cdots > b_n. \tag{3}$$

Hence, for $b_1 > \cdots > b_n$, the density function of $B_1$ equals

$$\begin{aligned} f(b_1, \dots, b_n) &= \lim_{x,y \to 0} \frac{\det(p_1(x_i, b_j))_{i,j=1}^n \cdot \det(p_1(b_i, y_j))_{i,j=1}^n}{\det(p_2(x_i, y_j))_{i,j=1}^n} \\ &= \frac{2^{n(n-1)/2}}{\pi^{n/2} \prod_{j=1}^{n-1} j!} \prod_{1 \le i < j \le n} |b_i - b_j|^2 \prod_{j-1}^{n} e^{-b_j^2}, \end{aligned} \tag{4}$$

which is the density function of the eigenvalues of $n \times n$ Hermitian matrices from the Gaussian unitary ensemble. Therefore, combined with the well-known results of random matrix theory,

$$\lim_{n \to \infty} \mathbb{P}((B_1^{(1)} - \sqrt{2n})\sqrt{2} n^{1/6} \le x) = F_{\mathrm{TW}}(x). \tag{5}$$

1.1.2. *Longest increasing subsequence and Plancherel measure on partitions.* The longest increasing subsequence problem can be formulated in the following manner. Denote by $S_n$ the symmetric group on $n$ symbols endowed with uniform measure. Given $\pi \in S_n$, a subsequence $\pi(i_1), \dots, \pi(i_r)$ is called an *increasing* subsequence if $i_1 < \cdots < i_r$ and $\pi(i_1) < \cdots < \pi(i_r)$. Denote by $\ell_n(\pi)$ the length of the longest increasing subsequence (this subsequence need not be unique). For applications of $\ell_n$ and activities around the asymptotic behavior of $\ell_n$, see, for example, [2, 5, 16]. In particular, in [5], the following limit theorem is:

$$\lim_{n \to \infty} \mathbb{P}\left( \frac{\ell_n(\pi) - 2\sqrt{n}}{n^{1/6}} < x \right) = F_{\mathrm{TW}}(x). \tag{6}$$

A closely related object is the uniform measure on the set of pairs of standard Young tableaux having the same shape (equivalently, the so-called Plancherel measure on the set of partitions). Given a partition of $n$, $\lambda = (\lambda_1, \dots, \lambda_r)$, where $\lambda_1 \ge \cdots \ge \lambda_r > 0$ and $\lambda_1 + \cdots + \lambda_r = n$, a standard Young tableaux of shape $\lambda$ consists of $r$ rows of boxes with distinct entries from $\{1, \dots, n\}$ such that the rows are left-justified, the $i$th row has $\lambda_i$ boxes and the entries are constrained to increase along rows and columns from left to right and top to bottom, respectively. These objects will be called *row-increasing Young tableaux* if the rows increase but the columns do not necessarily increase. The Robinson–Schensted bijection implies that the number of boxes in the top row of the pair of standard Young tableaux corresponding to $\pi \in S_n$ is equal to $\ell_n(\pi)$ [49]. Therefore, the distribution of $\ell_n$ is the same as the distribution of the number of boxes in the top row of the pair of standard Young tableaux having the same, shape chosen uniformly. This correspondence provides a representation of $\ell_n$ which is computable in terms of explicit formulae if the number of standard Young tableaux of a given shape is computable.

One way (among many) to compute the number of standard Young tableaux of shape $\lambda$ is by means of a nonintersecting path argument [35]. Let $N_t^1, \dots, N_t^r$ be independent rate-1 Poisson processes with initial conditions $N_0^i = 1 - i$ for $i = 1, 2, \dots, r$. Define $A_\lambda$ to be the event that $N_1^i = \lambda_i + (1 - i)$ for all $i = 1, 2, \dots, r$. For almost every element of $A_\lambda$ (the elements of $A_\lambda$ where no two jumps of these processes occur at the same time), there is a natural map to a row-increasing Young tableaux. The map is defined as follows. If $N^i$ jumps first, then place a 1 in the leftmost box in row $i$; if $N^j$ jumps second, then place a 2 in the first box of row $j$ if $j \neq i$ and a 2 in the second box of row $i$ if $j = i$. Continue in this fashion to produce a row-increasing Young tableaux of shape $\lambda$. It is not hard to show that this map induces the uniform probability measure [when properly normalized by $\mathbb{P}(A_\lambda)$] on the row-increasing Young tableaux. The subset $B_\lambda \subset A_\lambda$ which is mapped to the standard Young tableaux of shape $\lambda$ corresponds to the set of realizations whose paths do not intersect each other for all $t \in [0, 1]$. Since the mapping described induces uniform measure on the row-increasing Young tableaux of shape $\lambda$ and the standard Young tableaux correspond to nonintersecting path realizations, $B_\lambda$, the number of standard Young tableaux of shape $\lambda$, can be computed by evaluating

$$|\text{row-increasing Young tableaux of shape } \lambda| \frac{\mathbb{P}(B_\lambda)}{\mathbb{P}(A_\lambda)}. \tag{7}$$

The denominator of (7) is $e^{-r} \prod_{i=1}^r \frac{1}{\lambda_i!}$, by definition of Poisson processes and the independence of the $N^i$, while $|$ row-increasing Young tableaux of shape $\lambda| = \frac{n!}{\lambda_1! \cdots \lambda_r!}$ by elementary combinatorics. On the other hand, via the Karlin–McGregor formula [36],

$$\mathbb{P}(B_\lambda) = \det\left(\frac{e^{-1}}{(\lambda_i - i + j)!}\right)_{i,j=1}^r. \tag{8}$$

Hence, the number of standard Young tableaux of shape $\lambda$ is $n! \det(\frac{1}{(\lambda_i - i + j)!})_{i,j=1}^r$. In tandem with the RSK correspondence, this formula leads to an algebraic formula for the number of $\pi \in S_n$ for which $\ell_n(\pi) \leq m$. Moreover, a slight extension of this argument shows that result (6) can be stated in terms of the top curve of the nonintersecting Poisson processes if these processes were forced to return to their initial locations at time 2 by imposing that their dynamics between times 1 and 2 have negative rather than positive jumps. The asymptotic behavior of other curves can also be studied [6, 7, 13, 31, 43].

1.1.3. *Symmetric simple random walks and random rhombus tilings of a hexagon.* Consider $n$ symmetric simple (Bernoulli) random walks $S(m) = (S^{(1)}(m), \dots, S^{(n)}(m))$, conditioned not to intersect and such that $S(0) = (2(n-1), 2(n-2), \dots, 0) = S(2k)$. Any realization of such walks is in one-to-one correspondence with a rhombus tiling of a hexagon with side lengths $k, k, n, k, k, n$. Again, using the argument of Karlin and McGregor, the distribution of $S(k)$ can be expressed in terms of a determinant. This determinant was significantly simplified and was shown to be related to the so-called Hahn orthogonal polynomials by Johansson [33]. A further asymptotic analysis of the Hahn polynomials [8, 9] shows that as $n, k \to \infty$ such that $k = O(n)$, the top walk $S^{(1)}(k)$ converges to $F_{\mathrm{TW}}$ and the gap distribution "in bulk" converges to a discrete version of $F_{\mathrm{Sine}}$. A similar asymptotic result was also obtained for domino tilings of an Aztec diamond [33].

Certain polynuclear growth models, last passage percolation problems and a bus system problem [4, 34, 42, 46] have also been analyzed in depth using nonintersecting path techniques. In each of the cases described above, the random walks are very specific and the analysis relies heavily on their particular properties.

1.2. *Statement of theorems.* Let $k$ be a positive integer. Let

$$x_i = \frac{2i-k}{k}, \qquad i \in \{0, \dots, k\}. \tag{9}$$

Note that $x_i \in [-1, 1]$ for all $i$. Let $\{Y_l^j\}_{j=0,l=1}^{k,N_k}$ be a family of independent identically distributed random variables where $N_k$ is a positive integer. Assume that $\mathbb{E}Y_l^j = 0$ and $\mathrm{Var}(Y_l^j) = 1$. Further, assume that there exists $\lambda_0 > 0$ such that $\mathbb{E}(e^{\lambda Y_l^j}) < \infty$ for all $|\lambda| < \lambda_0$.

Define the random walk process $S(t) = (S_0(t), \dots, S_k(t))$ by

$$S_j(t) = x_j + \sqrt{\frac{2}{N_k}}\left(\sum_{i=1}^{\lfloor tN_k/2 \rfloor} Y_i^j + \left(\frac{tN_k}{2} - \left\lfloor \frac{tN_k}{2} \right\rfloor\right) Y^j_{\lfloor tN_k/2 \rfloor + 1}\right) \quad \text{for } t \in [0, 2], \tag{10}$$

which starts at $S_j(0) = x_j$. For $N_k$ equally spaced times, $S_j$ is given by

$$S_j\left(\frac{2}{N_k} l\right) = x_j + \sqrt{\frac{2}{N_k}}(Y_1^j + \dots + Y_l^j), \qquad l = 1, 2, \dots, N_k. \tag{11}$$

For $t$ between $\frac{2}{N_k} l$ and $\frac{2}{N_k}(l+1)$, $S_j(t)$ is simply defined by linear interpolation.

Let $(C([0,2]; \mathbb{R}^{k+1}), \mathcal{C})$ be the family of measurable spaces constructed from the continuous functions on $[0, 2]$ taking values in $\mathbb{R}^{k+1}$ equipped with their Borel sigma algebras (generated by the sup norm). Let $A_k, B_k \in \mathcal{C}$ be the events defined by

$$A_k = \{y_0(t) < \cdots < y_k(t) \text{ for } t \in [0,2]\}, \tag{12}$$

$$B_k = \{y_i(2) \in [x_i - h_k, x_i + h_k] \text{ for } i \in \{0,\ldots,k\}\}, \tag{13}$$

where $h_k > 0$. The results of this paper focus on the process $S(t)$ conditioned on the event $A_k \cap B_k$, where $h_k \ll \frac{2}{k}$. In other words, the particles never intersect and all particles essentially return to their original locations at the final time 2.

The main results of this paper state that under certain technical conditions on $h_k$ and $N_k$, as $k \to \infty$, the locations of the particles at the half time $(t = 1)$ behave statistically, after suitable scaling, like the eigenvalues of a large random Hermitian matrix from the Gaussian unitary ensemble. The conditions for $h_k$ and $N_k$ are that $\{h_k\}_{k>0}$ is a sequence of positive numbers and that $\{N_k\}_{k>0}$ is a sequence of positive integers satisfying

$$h_k \le (2k)^{-2k^2} \quad \text{and} \quad N_k \ge h_k^{-4(k+2)}. \tag{14}$$

Let $C_k, D_k \in \mathcal{C}$ be defined by

$$C_k = \left\{ y_k(1) \le \sqrt{2k} + \frac{\xi}{\sqrt{2}k^{1/6}} \right\}, \tag{15}$$

$$D_k = \left\{ y_i(1) \notin \left[ -\frac{\pi\eta}{\sqrt{2k}}, \frac{\pi\eta}{\sqrt{2k}} \right] \text{ for all } i \in \{0,\ldots,k\} \right\}, \tag{16}$$

where $\xi$ and $\eta > 0$ are fixed real numbers. The event $C_k$ is a constraint on the location of the rightmost particle and $D_k$ is the event that no particle is in a small neighborhood of the origin at time 1.

THEOREM 1 (Edge). *Let $\mathbb{P}_k$ be the probability measure induced on $(C([0,2]; \mathbb{R}^{k+1}), \mathcal{C})$ by the random walks $\{S(t) : t \in [0,2]\}$. Let $\{h_k\}_{k>0}$ and $\{N_k\}_{k>0}$ satisfy (14). Then*

$$\lim_{k\to\infty} \mathbb{P}_k(C_k | A_k \cap B_k) = F_{\mathrm{TW}}(\xi). \tag{17}$$

A similar theorem holds for the bulk.

THEOREM 2 (Bulk). *Let $\mathbb{P}_k$ be the probability measure induced on $(C([0,2]; \mathbb{R}^{k+1}), \mathcal{C})$ by the random walks $\{S(t) : t \in [0,2]\}$ and let $\{h_k\}_{k>0}$ and $\{N_k\}_{k>0}$ satisfy (14). Then*

$$\lim_{k\to\infty} \mathbb{P}_k(D_k | A_k \cap B_k) = F_{\mathrm{Sine}}(\eta). \tag{18}$$

The proofs have a two-step strategy. The first step is to show that under the conditions of the theorems, the process $S(t)$ is well approximated by nonintersecting Brownian bridge processes starting and ending at the same positions. This proof relies on the Komlos–Major–Tusnady (KMT) theorem. The second step is to compute the limiting distributions of the nonintersecting Brownian bridge processes and prove that these distributions are indeed $F_{\mathrm{TW}}$ or $F_{\mathrm{Sine}}$. This process is quite similar to the one discussed in Section 1.1.1, with the minor change that the Brownian bridge processes start and end at equally spaced locations, rather than at the same location. This change results in a Coulomb-gas density with the so-called Stieltjes–Wigert potential instead of the quadratic potential which appears in the GUE case. Such a nonintersecting Brownian bridge process was also considered in [24, 26] and the connection to the Stieltjes–Wigert potential was made in [26] in order to compute the partition function and the limiting density of states. However, the edge and bulk scaling limits of the system had not been worked out. This paper obtains the asymptotics of the orthogonal polynomials with respect to the Stieltjes–Wigert weight by using the Riemann–Hilbert method. As a consequence, the edge and bulk scaling limits are obtained.

The above theorems are proved under the condition that $N_k$ is large compared to $k+1$, the number of particles. This assumption ensures that the Brownian approximation of the random walks has a smaller effect than the nonintersecting condition. Although it is believed that the condition on $N_k$ is technical, it is not clear under which conditions on the random variables one has $N_k = O(k)$. For example, when $\{Y_l^j\}$ are Bernoulli, these results were proven even when $N_k = O(k)$ (see Section 1.1.3 above). This is because there is an integrability in this problem: the Karlin–McGregor argument applies directly because intersecting paths must be incident at some time. It is a challenge to find the optimal scaling such that a result of this nature holds for more general random variables. In other words, in what scaling regime does the exact Karlin–McGregor calculation essentially not matter?

This paper is organized as follows. The approximation by the Brownian bridge process is proved in Section 2. The asymptotic analysis of the Brownian bridge process (appearing in Section 2) is carried out in Section 3. Some other considerations such as finite-dimensional distributions and the modifications necessary to study random variables without finite moment generating functions are discussed in Section 4.

**2. Approximation by a Brownian bridge process.** Let $\{X_t\}_{t\geq 0}$ be the $\mathbb{R}^{k+1}$-valued stochastic process $X_t = (X_0(t), \ldots, X_k(t))$, where $X_j(t) = x_j + B_t^j$ for a family of $k+1$ independent standard Brownian motions $B_t^j$. The proof in this section relies on the Komlos–Major–Tusnady coupling of Brownian motions and random walks [38, 39] which can be stated in our setting as follows. With increments of the form $\{Y_l^j\}_{j=0,l=1}^{k,N_k}$ described in the Introduction, there exists a coupling such that

$$\mathbb{P}\left(\sup_{0\le l\le N_k}\left|S_i\left(\frac{2l}{N_k}\right)-X_i\left(\frac{2l}{N_k}\right)\right|>\frac{1}{\sqrt{N_k}}(c\log N_k+x)\right)\le e^{-ax} \tag{19}$$

for some fixed $a,c>0$ which depend only on the properties of the moment generating functions of the $\{Y_l^j\}_{j=0,l=1}^{k,N_k}$. Alternatively, (19) can be written as

$$\mathbb{P}\left(\sup_{0\le l\le N_k}\left|S_i\left(\frac{2l}{N_k}\right)-X_i\left(\frac{2l}{N_k}\right)\right|>\frac{c\log N_k}{\sqrt{N_k}}+y\right)\le e^{-ay\sqrt{N_k}}. \tag{20}$$

This fact immediately implies that

$$\begin{aligned}&\mathbb{P}\left(\sup_{0\le i\le k}\sup_{0\le l\le N_k}\left|S_i\left(\frac{2l}{N_k}\right)-X_i\left(\frac{2l}{N_k}\right)\right|>\frac{c\log N_k}{\sqrt{N_k}}+y\right)\\&\qquad\le(k+1)e^{-ay\sqrt{N_k}}.\end{aligned} \tag{21}$$

Let $\{S(t)\}_{t\in[0,2]}$ be the $(k+1)$-dimensional random walk process defined in the Introduction and let $\{X_t\}$ be the KMT-coupled $(k+1)$-dimensional Brownian process on the same probability spaces $(\Omega^{(k)},\mathcal{F}^{(k)},\mathbb{P}^{(k)})$. We can assume that the probability space which holds $S$ and $X$ is large enough to hold a third process $Z_t=(Z_0(t),\dots,Z_k(t))$, where the $Z_i(t)$ are standard Brownian bridge processes with initial and terminal conditions specified by $Z_i(0)=Z_i(2)=x_i$. Let $F_k^S,F_k^X,F_k^Z:C([0,2],\mathbb{R}^{k+1})\to\mathbb{R}$ be defined by

$$F_k^S(y)=\frac{\mathbb{1}_{A_k\cap B_k\cap C_k}(y)}{\mathbb{E}\mathbb{1}_{A_k\cap B_k}(S)}, \tag{22}$$

$$F_k^X(y)=\frac{\mathbb{1}_{A_k\cap B_k\cap C_k}(y)}{\mathbb{E}\mathbb{1}_{A_k\cap B_k}(X)}, \tag{23}$$

$$F_k^Z(y)=\frac{\mathbb{1}_{A_k\cap B_k\cap C_k}(y)}{\mathbb{E}\mathbb{1}_{A_k\cap B_k}(Z)}. \tag{24}$$

Let $G_k^S,G_k^X,G_k^Z:C([0,2],\mathbb{R}^{k+1})\to\mathbb{R}$ be defined by

$$G_k^S(y)=\frac{\mathbb{1}_{A_k\cap B_k\cap D_k}(y)}{\mathbb{E}\mathbb{1}_{A_k\cap B_k}(S)}, \tag{25}$$

$$G_k^X(y)=\frac{\mathbb{1}_{A_k\cap B_k\cap D_k}(y)}{\mathbb{E}\mathbb{1}_{A_k\cap B_k}(X)}, \tag{26}$$

$$G_k^Z(y)=\frac{\mathbb{1}_{A_k\cap B_k\cap D_k}(y)}{\mathbb{E}\mathbb{1}_{A_k\cap B_k}(Z)}. \tag{27}$$

Theorems 1 and 2 will be proven in two steps. The first step is to show that under the conditions given in the Introduction, the following holds.

PROPOSITION 1. *As $k \to \infty$, the random variables $F_k^S$, $F_k^Z$, $G_k^S$ and $G_k^Z$ satisfy*

$$\mathbb{E}(F_k^S(S) - F_k^Z(Z)) \to 0, \tag{28}$$

$$\mathbb{E}(G_k^S(S) - G_k^Z(Z)) \to 0. \tag{29}$$

As $\mathbb{E}(F_k^S(S)) = \mathbb{P}_k(C_k|A_k \cap B_k)$ and $\mathbb{E}(G_k^S(S)) = \mathbb{P}_k(D_k|A_k \cap B_k)$, it is enough to prove the following.

PROPOSITION 2.

$$\lim_{k\to\infty} \mathbb{E}F_k^Z(Z) = F_{\mathrm{TW}}(\xi), \tag{30}$$

$$\lim_{k\to\infty} \mathbb{E}G_k^Z(Z) = F_{\mathrm{Sine}}(\eta). \tag{31}$$

The proof of Proposition 2 is given in Section 3 below. The rest of this section focuses on the proof of Proposition 1. Proposition 1 is proved in two steps: first, $\mathbb{E}(F_k^S(S))$ is approximated by $\mathbb{E}(F_k^X(X))$ and second, $\mathbb{E}(F_k^X(X))$ is approximated by $\mathbb{E}(F_k^Z(Z))$. The proof of (29) is handled in a similar way.

Three preliminary lemmas are needed in order to prove Proposition 1. Recall from (9) that

$$x_i = \frac{2i-k}{k}, \qquad i = 0, \ldots, k. \tag{32}$$

LEMMA 1. *Let $a, c > 0$ be the constants in the KMT approximation (21). For any $\rho \geq \frac{3c\log N_k}{\sqrt{N_k}}$,*

$$\begin{aligned}\mathbb{E}|\mathbb{1}_{A_k\cap B_k}(S) - \mathbb{1}_{A_k\cap B_k}(X)| &\leq 2(k+1)e^{-(1/2)a\sqrt{N_k}\rho} \\ &\quad + \frac{32(k+1)}{\rho}\sqrt{\frac{N_k}{\pi}}e^{-\rho^2 N_k/64} + 8(2k+1)\rho,\end{aligned} \tag{33}$$

$$\begin{aligned}\mathbb{E}|\mathbb{1}_{B_k}(S) - \mathbb{1}_{B_k}(X)| &\leq (k+1)e^{-(1/2)a\sqrt{N_k}\rho} \\ &\quad + \frac{16(k+1)}{\rho}\sqrt{\frac{N_k}{\pi}}e^{-\rho^2 N_k/64} + 8(k+1)\rho,\end{aligned} \tag{34}$$

$$\begin{aligned}\mathbb{E}|\mathbb{1}_{C_k}(S) - \mathbb{1}_{C_k}(X)| &\leq (k+1)e^{-(1/2)a\sqrt{N_k}\rho} \\ &\quad + \frac{16(k+1)}{\rho}\sqrt{\frac{N_k}{\pi}}e^{-\rho^2 N_k/64} + 8(k+1)\rho,\end{aligned} \tag{35}$$

$$\mathbb{E}|\mathbb{1}_{D_k}(S)-\mathbb{1}_{D_k}(X)| \le (k+1)e^{-(1/2)a\sqrt{N_k}\rho} + \frac{16(k+1)}{\rho}\sqrt{\frac{N_k}{\pi}}e^{-\rho^2 N_k/64} + 8(k+1)\rho. \tag{36}$$

Proof. Note that

$$\begin{aligned}\mathbb{E}|\mathbb{1}_{A_k\cap B_k}(S)-\mathbb{1}_{A_k\cap B_k}(X)| &= \mathbb{E}|\mathbb{1}_{A_k}(S)\mathbb{1}_{B_k}(S)-\mathbb{1}_{A_K}(X)\mathbb{1}_{B_k}(X)| \\ &= \mathbb{E}|(\mathbb{1}_{A_k}(S)-\mathbb{1}_{A_k}(X))\mathbb{1}_{B_k}(S)+(\mathbb{1}_{B_k}(S)-\mathbb{1}_{B_k}(X))\mathbb{1}_{A_k}(X)| \\ &\le \mathbb{E}|\mathbb{1}_{A_k}(S)-\mathbb{1}_{A_k}(X)|+\mathbb{E}|\mathbb{1}_{B_k}(S)-\mathbb{1}_{B_k}(X)|.\end{aligned} \tag{37}$$

We first estimate $\mathbb{E}|\mathbb{1}_{A_k}(S)-\mathbb{1}_{A_k}(X)| = \mathbb{P}(\mathcal{E})$, where $\mathcal{E}=\{S\in A_k, X\notin A_k\}\cup\{S\notin A_k, X\in A_k\}$. Recall that $A_k=\{y_0(t)<\cdots<y_k(t) \text{ for } t\in[0,2]\}$. Let $\rho\ge\frac{3c\log N_k}{\sqrt{N}k}$, where $c$ is the KMT coupling constant. The event $\mathcal{E}$ can be expressed as the disjoint union of the three events $\mathcal{E}_1,\mathcal{E}_2,\mathcal{E}_3$. The first event $\mathcal{E}_1$ is the subset of $\mathcal{E}$ consisting of "bad" paths satisfying $\sup_{0\le i\le k}\sup_{t\in[0,2]}|S_i(t)-X_i(t)|>\frac{c\log N_k}{\sqrt{N_k}}+\rho$. The second event $\mathcal{E}_2$ is the subset of $\mathcal{E}\setminus\mathcal{E}_1$ consisting of paths satisfying $\min_{t\in[0,2]}(S_i(t)-S_{i-1}(t))<0$ for some $1\le i\le k$ while $X_0(t)<\cdots<X_k(t)$ for all $t\in[0,2]$. The third event $\mathcal{E}_3$ is the subset of $\mathcal{E}\setminus\mathcal{E}_1$ consisting of paths such that $\min_{t\in[0,2]}(X_i(t)-X_{i-1}(t))<0$ for some $1\le i\le k$ while $S_0(t)<\cdots<S_k(t)$ for all $t\in[0,2]$. In order to estimate $\mathbb{P}(\mathcal{E}_1)$, note that the KMT theorem couples random walks to Brownian motion at discrete times. Hence, even when $X$ and $S$ are close at discrete times, "bad paths" may occur if $X$ fluctuates too much in $(\frac{2}{N_k}l,\frac{2}{N_k}(l+1))$ for some $l$. (Note that $S$ is simply linearly interpolated for times not integral multiple of $\frac{2}{N_k}$.) Thus, from (21) and standard estimate for Brownian motions,

$$\begin{aligned}\mathbb{P}(\mathcal{E}_1) &\le \mathbb{P}\left(\sup_{0\le i\le k}\sup_{0\le l\le N_k}\left|S_i\left(\frac{2l}{N_k}\right)-X_i\left(\frac{2l}{N_k}\right)\right|>\frac{c\log N_k}{\sqrt{N_k}}+\frac{\rho}{2}\right) \\ &\quad+\mathbb{P}\Bigg(\left\{\sup_{0\le i\le k}\sup_{0\le l\le N_k}\left|S_i\left(\frac{2l}{N_k}\right)-X_i\left(\frac{2l}{N_k}\right)\right|\le\frac{c\log N_k}{\sqrt{N_k}}+\frac{\rho}{2}\right\} \\ &\qquad\cap\left\{\max_{s,t\in(2l/N_k,2(l+1)/N_k)}|X_i(t)-X_i(s)|>\frac{\rho}{2}\right. \\ &\qquad\qquad\left.\text{for some } 0\le i\le k \text{ and for some } 0\le l<N_k\right\}\Bigg) \\ &\le (k+1)e^{-(1/2)a\sqrt{N_k}\rho} \\ &\quad+(k+1)N_k\mathbb{P}\left(\max_{t,s\in[0,2/N_k]}|X_1(t)-X_1(s)|>\frac{\rho}{2}\right)\end{aligned} \tag{38}$$

$$\leq (k+1)e^{-(1/2)a\sqrt{N_k}\rho} + \frac{16(k+1)}{\rho}\sqrt{\frac{N_k}{\pi}}e^{-\rho^2 N_k/64}.$$

Note that this estimate does not use the fact that $\mathcal{E}_1$ is a subset of $\mathcal{E}$. For a path in the event $\mathcal{E}_2$, there exists $i \in \{1, 2, \dots, k\}$ such that $\min_{t\in[0,2]}(S_i(t) - S_{i-1}(t)) < 0$, but $X_{i-1}(t) < X_i(t)$ for all $t \in [0, 2]$ and $|S_j(t) - X_j(t)| \leq \frac{c \log N_k}{\sqrt{N_k}} + \rho$ all $t \in [0, 2]$ and $j \in \{0, 1, \dots, N_k\}$. Therefore, for a path in $\mathcal{E}_2$, $0 < \min_{t\in[0,2]}(X_i(t) - X_{i-1}(t)) < \frac{2c \log N_k}{\sqrt{N_k}} + 2\rho \leq 4\rho$. Thus, from a standard Brownian motion argument,

$$\begin{aligned} \mathbb{P}(\mathcal{E}_2) &\leq \mathbb{P}\Big(0 \leq \min_{t\in[0,2]}(X_i(t) - X_{i-1}(t)) < 4\rho \text{ for some } 1 \leq i \leq k\Big) \\ &\leq k\mathbb{P}\Big(0 \leq \min_{t\in[0,2]}(X_1(t) - X_0(t)) < 4\rho\Big) \leq 4k\rho. \end{aligned} \tag{39}$$

A similar argument yields that

$$\begin{aligned} \mathbb{P}(\mathcal{E}_3) &\leq \mathbb{P}\Big(-4\rho < \min_{t\in[0,2]}(X_i(t) - X_{i-1}(t)) < 0 \text{ for some } 1 \leq i \leq k\Big) \\ &\leq k\mathbb{P}\Big(-4\rho < \min_{t\in[0,2]}(X_1(t) - X_0(t)) < 0\Big) \leq 4k\rho. \end{aligned} \tag{40}$$

Therefore,

$$\begin{aligned} &\mathbb{E}|\mathbb{1}_{A_k}(S) - \mathbb{1}_{A_k}(X)| \\ &\quad = \mathbb{P}(\mathcal{E}_1) + \mathbb{P}(\mathcal{E}_2) + \mathbb{P}(\mathcal{E}_3) \\ &\quad \leq (k+1)e^{-(1/2)a\sqrt{N_k}\rho} + \frac{16(k+1)}{\rho}\sqrt{\frac{N_k}{\pi}}e^{-\rho^2 N_k/64} + 8k\rho. \end{aligned} \tag{41}$$

We now estimate $\mathbb{E}|\mathbb{1}_{B_k}(S) - \mathbb{1}_{B_k}(X)| = \mathbb{P}(\mathcal{F})$, where $\mathcal{F} = \{S \in B_k, X \notin B_k\} \cup \{S \notin B_k, X \in B_k\}$. As before, we express $\mathcal{F}$ as $\mathcal{F}_1 \cup \mathcal{F}_2 \cup \mathcal{F}_3$, a disjoint union. The first event, $\mathcal{F}_1$, is the subset of $\mathcal{F}$ consisting of the same bad paths as in $\mathcal{E}_1$. The event $\mathcal{F}_2$ is the intersection of $\mathcal{F} \setminus \mathcal{F}_1$ and $\{S \in B_k, X \notin B_k\}$ and the event $\mathcal{F}_3$ is the intersection of $\mathcal{F} \setminus \mathcal{F}_1$ and $\{S \notin B_k, X \in B_k\}$. The argument for $\mathcal{E}_1$ implies that the same bound (38) applies to $\mathbb{P}(\mathcal{F}_1)$. For a path in $\mathcal{F}_2$, there exists $i \in \{0, 1, \dots, k\}$ such that $X_i(2) \notin [x_i - h_k, x_i + h_k]$. But as $S_i(2) \in [x_i - h_k, x_i + h_k]$ and $|S_i(2) - X_i(2)| \leq \frac{2 \log N_k}{\sqrt{N_k}} + \rho \leq 2\rho$, we find that $X_i(2) \in (x_i + h_k, x_i + h_k + 2\rho]$ or $X_i(2) \in [x_i - h_k - 2\rho, x_i - h_k)$. Therefore,

$$\mathbb{P}(\mathcal{F}_2) \leq (k+1)\mathbb{P}(X_0(2) \in [-2\rho, 2\rho]) \leq 4(k+1)\rho. \tag{42}$$

A similar argument yields the same bound for $\mathbb{P}(\mathcal{F}_3)$. Therefore, (34) is proved, as is (33), by using (37) and (41). An almost identical argument proves (35) and (36). $\square$

Denote by $p_t(a,b) = \frac{1}{\sqrt{2\pi t}} e^{-(a-b)^2/(2t)}$ the standard heat kernel in one dimension. The theorem of Karlin and McGregor [36] for nonintersecting Brownian motions implies that the joint probability density function $f_t(y_0,\dots,y_k)$ of $(k+1)$-dimensional Brownian motion $X(t)$ at time $t$ satisfying $X_0(s) < X_1(s) < \cdots < X_k(s)$ for $s \in [0,t]$ is equal to

$$f_t(y_0,\dots,y_k) = \det(p_t(x_i,y_j))_{i,j=0}^k, \tag{43}$$

where $x_i = X_i(0)$. The following lemma establishes a lower bound for this density when $y_i = x_i$ for all $i$.

LEMMA 2. *For $t > 0$,*

$$\det(p_t(x_i,x_j))_{i,j=0}^k \ge \frac{1}{(2\pi t)^{(k+1)/2}} e^{-2(k+1)(k+2)/(3tk)} \left(\frac{2}{\sqrt{t}k}\right)^{k(k+1)}. \tag{44}$$

*In particular, for all sufficiently large $k$,*

$$\det(p_2(x_i,x_j))_{i,j=0}^k \ge k^{-k^2}. \tag{45}$$

PROOF. As $x_i = \frac{2i-k}{k}$, we have

$$\begin{aligned} \det(p_t(x_i,x_j))_{i,j=0}^k &= \det\left(\frac{1}{\sqrt{2\pi t}} e^{-1/(2t)(x_i-x_j)^2}\right)_{i,j=0}^k \\ &= \frac{e^{-2\sum_{j=0}^k j^2}}{(2\pi t)^{(k+1)/2}} \det(e^{2ij/(tk^2)})_{i,j=0}^k. \end{aligned} \tag{46}$$

It is an exercise to show that for $k \ge 1$,

$$\det(e^{2ij/(tk^2)})_{i,j=0}^k = \left[\prod_{l=1}^k \delta^{l(l-1)/2}\right]\left[\prod_{j=1}^k (\delta^j - 1)^{k+1-j}\right], \tag{47}$$

where $\delta = e^{4/(tk^2)}$. Using (47) and the fact that $\delta - 1 > \frac{4}{tk^2} > 0$, we obtain

$$\begin{aligned} &\det(p_t(x_i,x_j))_{i,j=0}^k \\ &\quad = \frac{1}{(2\pi t)^{(k+1)/2}} e^{-2(k+1)(k+2)/(3tk)} \prod_{j=1}^k (\delta^j - 1)^{k+1-j} \\ &\quad \ge \frac{1}{(2\pi t)^{(k+1)/2}} e^{-2(k+1)(k+2)/(3tk)} (\delta - 1)^{k(k+1)/2} \\ &\quad \ge \frac{1}{(2\pi t)^{(k+1)/2}} e^{-2(k+1)(k+2)/(3tk)} \left(\frac{4}{tk^2}\right)^{k(k+1)/2}. \end{aligned} \tag{48}$$

This completes the proof of Lemma 2. □

The following lemma will be used to control the difference between a conditioned version of the process $X$ and the process $Z$.

LEMMA 3. *If* $h_k \le (2k)^{-2k^2}$, *then for sufficiently large* $k$,

$$\left|\frac{\det(p_1(y_i,x_j))}{\det(p_2(x_i,x_j))} - \frac{\int_{-h_k}^{h_k}\cdots\int_{-h_k}^{h_k}\det(p_1(y_i,x_j+s_j))\,ds_0\cdots ds_k}{\int_{-h_k}^{h_k}\cdots\int_{-h_k}^{h_k}\det(p_2(x_i,x_j+s_j))\,ds_0\cdots ds_k}\right| \le \frac{1}{k}, \tag{49}$$

*uniformly in* $(y_0,\dots,y_k)\in\mathbb{R}^{k+1}$.

PROOF. The conclusion of this lemma is a consequence of several elementary determinant estimates. First, note that if $A=(a_{ij})_{i,j=0}^k$ is a $(k+1)\times(k+1)$ matrix with entries $|a_{ij}|\le 1$, then for the matrix $I^{ij}$ given by $(I^{ij})_{mn}=\delta_{im}\delta_{jn}$, we have

$$|\det(A)-\det(A+\varepsilon I^{ij})| \le \varepsilon k!. \tag{50}$$

Using a Lipschitz estimate for the Gaussian density, equation (50) implies that for any $t\ge\frac{1}{\sqrt{2\pi e}}$, any $h>0$ and any $(a_0,\dots,a_k),(b_0,\dots,b_k)\in\mathbb{R}^{k+1}$,

$$\begin{aligned}&\left|\det(p_t(a_i,b_j)) - \frac{1}{(2h)^{k+1}}\int_{-h}^{h}\cdots\int_{-h}^{h}\det(p_t(a_i+s_i,b_j))\,ds_0\cdots ds_k\right|\\ &\qquad \le 2h(k+1)^2k!.\end{aligned} \tag{51}$$

A simple algebraic manipulation now yields that the left-hand side of (49) equals

$$\left|\frac{\det(p_1(y_i,x_j))}{\det(p_2(x_i,x_j))}\cdot\frac{\mathcal{Q}_2}{\det(p_2(x_i,x_j))+\mathcal{Q}_2} + \frac{\mathcal{Q}_1}{\det(p_2(x_i,x_j))+\mathcal{Q}_2}\right|, \tag{52}$$

where

$$\begin{aligned}\mathcal{Q}_1 = \frac{1}{(2h_k)^{k+1}}\int_{[-h_k,h_k]^{k+1}}&[\det(p_1(y_i,x_j))\\ &-\det(p_1(y_i,x_j+s_j))]\,ds_0\cdots ds_k,\end{aligned} \tag{53}$$

$$\begin{aligned}\mathcal{Q}_2 = \frac{1}{(2h_k)^{k+1}}\int_{[-h_k,h_k]^{k+1}}&[\det(p_2(x_i,x_j+s_j))\\ &-\det(p_2(x_i,x_j))]\,ds_0\cdots ds_k.\end{aligned} \tag{54}$$

Using the estimates (45) and (51), we obtain

$$\det(p_2(x_i,x_j))+\mathcal{Q}_2 \ge k^{-k^2} - 2h_k(k+1)^2k! \ge \tfrac12 k^{-k^2} \tag{55}$$

for sufficiently large $k$. Hence, again using (51),

$$\left|\frac{\mathcal{Q}_1}{\det(p_2(x_i,x_j))+\mathcal{Q}_2}\right| \le 2h_k(k+1)^2k!k^{k^2}. \tag{56}$$

On the other hand, as $\det(p_1(x_i, y_j))$ is the density function for $(y_0, \ldots, y_k) \in \mathbb{R}^{k+1}_>$, where $\mathbb{R}^{k+1}_> = \{(y_0, \ldots, y_k) \in \mathbb{R}^{k+1} : y_0 < \cdots < y_k\}$ corresponding to the probability of $k+1$ Brownian motions starting from $(x_0, \ldots, x_k)$ and ending at $(y_0, \ldots, y_k)$ at time 1 without having intersected, it is clearly less than the same type of probability density function when a nonintersection condition is not imposed. Therefore,

$$\det(p_1(x_i, y_j)) \le \prod_{i=0}^{k} \frac{1}{\sqrt{2\pi}} e^{-(1/2)(x_i - y_i)^2} \le 1 \tag{57}$$

and hence

$$\left| \frac{\det(p_1(y_i, x_j))}{\det(p_2(x_i, x_j))} \cdot \frac{\mathcal{Q}_2}{\det(p_2(x_i, x_j)) + \mathcal{Q}_2} \right| \le 2h_k (k+1)^2 k! k^{2k^2}. \tag{58}$$

Since $h_k$ is assumed to be less than or equal to $(2k)^{-2k^2}$, (49) follows. □

Proof of Proposition 1. Two estimates will be needed. Note that

$$|\mathbb{E}(F_k^S(S) - F_k^Z(Z))| \le \mathbb{E}|F_k^S(S) - F_k^X(X)| + |\mathbb{E}(F_k^X(X) - F_k^Z(Z))|. \tag{59}$$

The first term on the right-hand side of (59) is estimated as follows:

$$\begin{aligned}
&\mathbb{E}|F_k^S(S) - F_k^X(X)| \\
&\qquad = \mathbb{E}\left| \frac{\mathbb{1}_{A_k \cap B_k \cap C_k}(S)}{\mathbb{E}\mathbb{1} A_k \cap B_k(S)} - \frac{\mathbb{1}_{A_k \cap B_k \cap C_k}(X)}{\mathbb{E}\mathbb{1}_{A_k \cap B_k}(X)} \right| \\
&\qquad = \mathbb{E}|(\mathbb{1}_{A_k \cap B_k \cap C_k}(S)(\mathbb{E}\mathbb{1}_{A_k \cap B_k}(X) - \mathbb{E}\mathbb{1}_{A_k \cap B_k}(S)) \\
&\quad + (\mathbb{1}_{A_k \cap B_k \cap C_k}(S) - \mathbb{1}_{A_k \cap B_k \cap C_k}(X))\mathbb{E}\mathbb{1}_{A_k \cap B_k}(S)) \\
&\qquad\qquad \times (\mathbb{E}\mathbb{1}_{A_k \cap B_k}(S)\mathbb{E}\mathbb{1}_{A_k \cap B_k}(X))^{-1}| \\
&\le \frac{|\mathbb{E}\mathbb{1}_{A_k \cap B_k}(X) - \mathbb{E}\mathbb{1}_{A_k \cap B_k}(S)|}{\mathbb{E}\mathbb{1}_{A_k \cap B_k}(X)} + \mathbb{E}\left| \frac{(\mathbb{1}_{A_k \cap B_k \cap C_k}(S) - \mathbb{1}_{A_k \cap B_k \cap C_k}(X))}{\mathbb{E}\mathbb{1}_{A_k \cap B_k}(X)} \right| \\
&\le \frac{2\mathbb{E}|\mathbb{1}_{A_k \cap B_k}(S) - \mathbb{1}_{A_k \cap B_k}(X)| + \mathbb{E}|\mathbb{1}_{C_k}(S) - \mathbb{1}_{C_k}(X)|}{\mathbb{E}\mathbb{1}_{A_k \cap B_k}(X)}.
\end{aligned} \tag{60}$$

By setting $\rho = N_k^{-1/4}$ in Lemma 1, for sufficiently large $k$, it is easy to check that

$$\begin{aligned}
\mathbb{E}|\mathbb{1}_{A_k \cap B_k}(S) - \mathbb{1}_{A_k \cap B_k}(X)| &\le \frac{20k}{N_k^{1/4}}, \\
\mathbb{E}|\mathbb{1}_{C_k}(S) - \mathbb{1}_{C_k}(X)| &\le \frac{20k}{N_k^{1/4}}.
\end{aligned} \tag{61}$$

On the other hand, by using (45) and the argument leading to (55), for sufficiently large $k$,

$$
\begin{aligned}
&\mathbb{E}\mathbb{1}_{A_k\cap B_k}(X)\\
&\quad= (2h_k)^{k+1}\det(p_2(x_i,x_j))\\
&\qquad+ (\mathbb{E}\mathbb{1}_{A_k\cap B_k}(X) - (2h_k)^{k+1}\det(p_2(x_i,x_j)))\\
&\quad= (2h_k)^{k+1}\det(p_2(x_i,x_j))\\
&\qquad+ \int_{[-h_k,h_k]^{k+1}} (\det(p_2(x_i,x_j+s_j)) - \det(p_2(x_i,x_j)))\,ds_0\cdots ds_k\\
&\quad\geq \frac{(2h_k)^{k+1}}{2k^{k^2}}.
\end{aligned}
\tag{62}
$$

Hence, from (61), for sufficiently large $k$,

$$
\mathbb{E}|F_k^S(S) - F_k^X(X)| \leq \frac{120k^{k^2+1}}{(2h_k)^{k+1}N_k^{1/4}} \to 0
\tag{63}
$$

as $k\to\infty$. For the second term of (59), note that the Karlin–McGregor formula for nonintersecting Brownian motions implies that [cf. (43) above] the density function of the nonintersecting Brownian bridge process $Z$ evaluated at time 1 is equal to

$$
f(y_0,\ldots,y_k) = \frac{\det(p_1(x_i,y_j))_{i,j=0}^k \det(p_1(y_i,x_j))_{i,j=0}^k}{\det(p_2(x_i,x_j))_{i,j=0}^k}.
\tag{64}
$$

Similarly, the density of the nonintersecting Brownian motion $X$ evaluated at time $t$ is equal to

$$
\begin{aligned}
&f(y_0,\ldots,y_k)\\
&\quad= \frac{\int_{[-h_k,h_k]^{k+1}} \det(p_1(x_i,y_j))_{i,j=0}^k \det(p_1(y_i,x_j+s_j))_{i,j=0}^k\,ds_0\cdots ds_k}{\int_{[-h_k,h_k]^{k+1}} \det(p_2(x_i,x_j+s_j))\,ds_0\cdots ds_k}.
\end{aligned}
\tag{65}
$$

Therefore,

$$
\begin{aligned}
&|\mathbb{E}(F_k^X(X) - F_k^Z(Z))|\\
&\quad= \left|\mathbb{E}\left(\frac{\mathbb{1}_{A_k\cap B_k\cap C_k}(X)}{\mathbb{E}\mathbb{1}_{A_k\cap B_k}(X)} - \frac{\mathbb{1}_{A_k\cap B_k\cap C_k}(Z)}{\mathbb{E}\mathbb{1}_{A_k\cap B_k}(Z)}\right)\right|\\
&\quad\leq \int_{\mathbb{R}_>^{k+1}} \left|\left(\int_{[-h_k,h_k]^{k+1}} \det(p_1(x_i,y_j))\right.\right.\\
&\qquad\qquad\left.\times \det(p_1(y_i,x_j+s_j))\,ds_0\cdots ds_k\right)
\end{aligned}
\tag{66}
$$

$$\times\left(\int_{[-h_k,h_k]^{k+1}}\det(p_2(x_i,x_j+s_j))\,ds_0\cdots ds_k\right)^{-1}$$
$$-\frac{\det(p_1(x_i,y_j))\det(p_1(y_i,x_j))}{\det(p_2(x_i,x_j))}\Bigg|\,dy_0\cdots dy_k.$$

By using Lemma 3 and (57), this implies that

$$\begin{aligned}|\mathbb{E}(F_k^X(X)-F_k^Z(Z))| &\le \frac{1}{k}\int_{\mathbb{R}_>^{k+1}}|\det(p_1(x_i,y_j))|\,dy_0\cdots dy_k\\ &\le \frac{1}{k}\int_{\mathbb{R}_>^{k+1}}\left[\prod_{i=0}^{k}\frac{1}{\sqrt{2\pi}}e^{-(x_i-y_i)^2/2}\right]dy_0\cdots dy_k\le\frac{1}{k}.\end{aligned}\tag{67}$$

The proof of (29) is exactly the same. This completes the proof of Proposition 1. $\square$

**3. Asymptotics of a Brownian bridge process.** We prove Proposition 2 in this section. Together with the results of Section 2, this completes the proofs of Theorem 1 and Theorem 2.

From the density formula of Karlin and McGregor for a nonintersecting Brownian bridge processes [36] [cf. (43)],

$$\begin{aligned}\mathbb{E}(F_k^Z(Z)) &= \frac{1}{\det(p_2(x_i,x_j))_{i,j=0}^k}\\ &\quad\times\int_{\mathbb{R}_>^{k+1}}[\det(p_1(x_i,y_j))_{i,j=0}^k]^2\prod_{j=0}^{k}(1-\mathcal{H}_1(y_j))\,dy_j,\end{aligned}\tag{68}$$

where $x_i=\frac{2i-k}{k}$ and

$$\mathcal{H}_1(y)=\mathbb{1}_{(\sqrt{2k}+\xi/(\sqrt{2}k^{1/6}),\infty)}(y).\tag{69}$$

Also,

$$\begin{aligned}\mathbb{E}(G_k^Z(Z)) &= \frac{1}{\det(p_2(x_i,x_j))_{i,j=0}^k}\\ &\quad\times\int_{\mathbb{R}_>^{k+1}}[\det(p_1(x_i,y_j))_{i,j=0}^k]^2\prod_{j=0}^{k}(1-\mathcal{H}_2(y_j))\,dy_j,\end{aligned}\tag{70}$$

where

$$\mathcal{H}_2(y)=\mathbb{1}_{[-\eta/\sqrt{k+1},\eta/\sqrt{k+1}]}(y).\tag{71}$$

We need the limit of (68) and (70) as $k\to\infty$.

In the discussion below, $\mathcal{H}(y)$ denotes either $\mathcal{H}_1$ or $\mathcal{H}_2$. Indeed, the algebra below works for arbitrary bounded functions $\mathcal{H}(y)$. Using the formula for $p_t$ and the definition of $x_i$, an elementary algebraic manipulation using Vandermode determinants yields that (68) and (70) are equal to

$$C'_k \cdot \int_{\mathbb{R}^{k+1}} \prod_{0\le i<j\le k} (e^{2y_j/k} - e^{2y_i/k})^2 \prod_{j=0}^{k} (1-\mathcal{H}(y_j)) e^{-y_j^2-2y_j}\, dy_j, \tag{72}$$

where $C'_k$ is the normalization constant so that (72) becomes 1 when $\mathcal{H}(y) \equiv 0$:

$$C'_k = \frac{e^{-(k+1)(k+2)/(3k)}}{\det(p_2(x_i,x_j))_{i,j=0}^{k}(k+1)!(2\pi)^{k+1}}. \tag{73}$$

Note that the integration domain is changed to $\mathbb{R}^{k+1}$ by using the symmetry of the integrand. A similar calculation for the case when the assumption that the particles start at equally-spaced locations and arrive at $(y_0,\ldots,y_k)$ at time 1 (without any assumption regarding what happens after time 1) can be found in equation (4.7) of [24]. Introducing the change of variables $y_j = \frac{k}{2}\log u_j - 1$, (72) becomes

$$C_k \cdot \int_{\mathbb{R}_+^{k+1}} \prod_{0\le i<j\le k} (u_j-u_i)^2 \prod_{j=0}^{k} (1-\hat{\mathcal{H}}(u_j)) \frac{1}{u_j} e^{-(k^2/4)(\log u_j)^2}\, du_j, \tag{74}$$

where

$$\hat{\mathcal{H}}(u) = \mathcal{H}\Big(\frac{k}{2}\log u - 1\Big) \tag{75}$$

and the normalization constant is given by

$$C_k = \frac{k^{k+1} e^{-(k+1)(k+2)/(3k)}}{\det(p_2(x_i,x_j))_{i,j=0}^{k}(k+1)!(2\pi)^{k+1}}. \tag{76}$$

This is the standard $\beta = 2$ ensemble in random matrix theory on the real half-line $\mathbb{R}_+$ with the weight

$$w(u) = \frac{1}{u} e^{-(k^2/4)(\log u)^2} = e^{-(k^2/4)(\log u)^2 - \log u}. \tag{77}$$

Note that $w(u) = o(u^{-m})$ for any $m \ge 0$ as $u \to +\infty$, and $w(u) = o(u^m)$ for any $m \ge 0$ as $u \downarrow 0$.

With the change of variables $u = e^{-2/k^2}x$, (77) becomes

$$w(u)\, du = c \cdot e^{-(k^2/4)(\log x)^2}\, dx, \qquad c = e^{1/k^2}. \tag{78}$$

This is, up to a constant, the Stieltjes–Wigert weight, which is defined as

$$\pi^{-1/2} k e^{-k^2(\log x)^2} \tag{79}$$

(see, e.g., Section 2.7 of [51] or Section 3.27 of [37]). The moments for the Stieltjes–Wigert weight is an example of an indeterminate moment problem; hence, there are several weights that have the same moments as the weight (79). Another interesting feature of the Stieltjes–Wigert weight (79) is that the corresponding orthogonal polynomials (called *Stieltjes–Wigert polynomials*) are examples of so-called $q$-polynomials with $q = e^{-1/(2k^2)}$ (see, e.g., Section 3.27 of [37]). The connection between the nonintersecting Brownian bridge process $Z$ and the Stieltjes–Wigert weight was first observed in [26]; the Stieltjes–Wigert weight also appears in [28].

Various $\beta = 2$ matrix ensembles of the form (74) (on both the real line and subsets of the real line) have been analyzed asymptotically and it has been proven that the local statistics of the "eigenvalues" (or the particles $u_0, \dots, u_k$) are generically independent of the potential $w$. For example, such "universality" is proved when $w(x) = e^{-(k+1)V(x)}$ for an analytic weight $V$ on $\mathbb{R}$ or $\mathbb{R}_+$ satisfying certain growth conditions as $x \to \pm\infty$ (and as $x \to 0$ for weights on $\mathbb{R}_+$) (e.g., [11, 21, 40, 44]) and when $w(x) = e^{-Q(x)}$, where $Q(x)$ is a polynomial (e.g., [20]). However, the asymptotic analysis of the ensemble with the weight given in (77) above does not seem to appear in the literature. It is well known [see (80) below] that the asymptotics of $\beta = 2$ ensembles amount to the asymptotic analysis of the corresponding orthogonal polynomials. For our case, we need the asymptotics of the orthogonal polynomials of degree $k$ and $k + 1$ with respect to the weight (77) as $k \to \infty$; note that the weight also varies as $k$ increases. The asymptotics of Stieltjes–Wigert polynomials were recently studied in [29] and [55], but in different asymptotic regimes: the degree goes to infinity while the weight is fixed. Therefore, the analysis of this section seems to yield new results for asymptotics of Stieltjes–Wigert polynomials. Nevertheless, the asymptotic analysis of the orthogonal polynomials and the ensemble (74) with varying weight (77) can be done in a very similar way to the analysis in [20, 21] using the Deift–Zhou steepest-descent method for related Riemann–Hilbert problems (RHP's), which is now one of standard tools for asymptotic analysis of orthogonal polynomials. We note that [55] also used the Deift–Zhou method (for a different asymptotic regime) and our analysis has some overlap with the analysis of [55]. In this section, we present only a sketch of the analysis.

It is a standard result in random matrix theory (see, e.g., [41, 52]) that (74) equals

$$\det(1 - \mathbf{K}_k \hat{\mathcal{H}}), \tag{80}$$

where

$$\mathbf{K}_k(x, y) = \sqrt{w(x)w(y)} \frac{\gamma_k}{\gamma_{k+1}} \frac{p_{k+1}(x)p_k(y) - p_k(x)p_{k+1}(y)}{x - y} \tag{81}$$

is the Christoffel–Darboux kernel in which $p_n(x) = \gamma_n x^n + \cdots$ is the $n$th orthonormal polynomial with respect to $w$. Hence,

$$\begin{aligned}\mathbb{E}(F_k^Z(Z)) &= \det(1 - \mathbf{K}_k \hat{H}_1),\\ \mathbb{E}(G_k^Z(Z)) &= \det(1 - \mathbf{K}_k \hat{H}_2).\end{aligned} \tag{82}$$

Let $\mathbf{Y}(z)$ be the solution to the following Riemann–Hilbert problem: $\mathbf{Y}(x)$ is the $2\times 2$-matrix-valued function on $\mathbb{C}\setminus\overline{R}_+$ satisfying

- $\mathbf{Y}(z)$ is analytic for $z \in \mathbb{C}\setminus\overline{\mathbb{R}}_+$, $\mathbf{Y}_\pm(z) = \lim_{\varepsilon\downarrow 0}\mathbf{Y}(z \pm i\varepsilon)$ is continuous for $z \in \mathbb{R}_+$ and $\mathbf{Y}(z)$ is bounded as $z \to 0$;
- for $z \in \mathbb{R}_+$,

$$\mathbf{Y}_+(z) = \mathbf{Y}_-(z)\begin{pmatrix} 1 & w(z)\\ 0 & 1\end{pmatrix}; \tag{83}$$

- $\mathbf{Y}(z) z^{-(k+1)\sigma_3} = (\mathbf{I} + O(z^{-1}))$ uniformly as $z \to \infty$ such that $z \in \mathbb{C}\setminus\overline{\mathbb{R}}_+$, where $\sigma_3 = \left(\begin{smallmatrix} 1 & 0\\ 0 & -1\end{smallmatrix}\right)$.

There is a unique solution $\mathbf{Y}$ to this RHP and, in particular, the (11) and (21) entries of $\mathbf{Y}(z)$ are given by $\mathbf{Y}_{11}(z) = \gamma_{k+1}^{-1} p_{k+1}(z)$ and $\mathbf{Y}_{21}(z) = -2\pi i \gamma_k p_k(z)$ [25]. Note that the existence of $\mathbf{Y}$ under the condition that $\mathbf{Y}(z)$ is bounded as $z \to 0$ (rather than, e.g., that $\mathbf{Y}_{12}(z) = O(z^{-1})$, as in, say, [54]) is due to the fact that $w(x) \to 0$ faster than any polynomial as $x \to 0$. Thus, the Christoffel–Darboux kernel can be written as, by using the fact that $\det \mathbf{Y}(z) = 1$,

$$\mathbf{K}_k(x,y) = \sqrt{w(x)w(y)}\,\frac{1}{2\pi i(x-y)}\begin{pmatrix} 0 & 1\end{pmatrix}\mathbf{Y}^{-1}(y)\mathbf{Y}(x)\begin{pmatrix} 1\\ 0\end{pmatrix}. \tag{84}$$

One of the main ingredient in analyzing the RHP for orthogonal polynomials asymptotically is the so-called equilibrium measure and the corresponding "$g$-function." Let $\psi(x)\,dx$ be a measure on $\mathbb{R}_+ = \operatorname{supp}(w)$ with total mass

$$\int \psi(x)\,dx = k+1. \tag{85}$$

Define the "$G$-function"

$$G(z) = \int \log(z-x)\psi(x)\,dx, \qquad z \in \mathbb{C}\setminus\overline{\mathbb{R}}_+, \tag{86}$$

where log represents the log function on the standard branch so that $\log u = \log|u| + i\arg(u)$, where $|\arg(u)| < \pi$. It is customary to define $\psi$ to be the probability measure and define the $g$-function as in (86) [hence, $G(z) = (k+1)g(z)$], but in this paper, we use the above convention since it simplifies some formulas below. Note that

$$G_+(x) + G_-(x) = 2\int \log|x-y|\psi(y)\,dy, \qquad x \in \mathbb{R}_+. \tag{87}$$

We look for $G$ satisfying the following two conditions: there exists a constant $\ell$ such that

- $G_+(x)+G_-(x)+\log(w(x))-\ell=0$ for $x\in \mathrm{supp}(\psi)$,
- $G_+(x)+G_-(x)+\log(w(x))-\ell<0$ for $x\in \mathbb{R}_+\setminus \mathrm{supp}(\psi)$.

For such $G$, the measure $\psi$ is called the *equilibrium measure*.

Using the standard procedure to solve this variational problem (see, e.g., [19, 47]), one can compute the equilibrium measure for the weight (77).

LEMMA 4. *For the weight (77), the support of the equilibrium measure is* $[\mathbf{a},\mathbf{b}]$*, where*

$$\sqrt{\mathbf{a}}=e^{(2k+1)/k^2}-\sqrt{e^{(4k+2)/k^2}-e^{2/k}},$$
$$\sqrt{\mathbf{b}}=e^{(2k+1)/k^2}+\sqrt{e^{(4k+2)/k^2}-e^{2/k}}. \tag{88}$$

*The equilibrium measure is, for* $x\in[\mathbf{a},\mathbf{b}]$,

$$\psi(x)=\frac{1}{2\pi}\sqrt{(\mathbf{b}-x)(x-\mathbf{a})}h(x),\qquad h(z)=\frac{1}{2\pi i}\oint_C\frac{-(\log(w(s)))'}{(s-z)R(s)}ds, \tag{89}$$

*where* $R(z)=((z-\mathbf{a})(z-\mathbf{b}))^{1/2}$ *denotes the principal branch of the square root function and the simple closed contour* $C$ *contains* $z$ *and* $[\mathbf{a},\mathbf{b}]$*, inside does not touch* $(-\infty,0]$ *and is oriented counterclockwise. A residue calculation yields that*

$$\psi(x)=\frac{k^2}{2\pi x}\arctan\left(\frac{\sqrt{(\mathbf{b}-x)(x-\mathbf{a})}}{\sqrt{\mathbf{ab}}+x}\right),\qquad x\in[\mathbf{a},\mathbf{b}]. \tag{90}$$

We remark that $\mathbf{a}$ and $\mathbf{b}$ are sometimes called the Mhaskar–Rakhmanov–Saff numbers. The above $\mathbf{a}$ and $\mathbf{b}$ are obtained in [55]: with $\alpha_n$ and $\beta_n$ as in (2.2) and (2.3) of [55], we have

$$\mathbf{a}=(e^{-1/(2k^2)}\alpha_n)|_{k\mapsto k/2,n=k+1},\qquad \mathbf{b}=(e^{-1/(2k^2)}\beta_n)|_{k\mapsto k/2,n=k+1}. \tag{91}$$

Given this $\psi$, $G(z)$ is defined as in (86) and $\ell$ is defined as $\ell=2G(\mathbf{b})-\log(w(\mathbf{b}))=2G(\mathbf{a})-\log(w(\mathbf{a}))$. The function $h(z)$ in (89) is analytic in $z\in\mathbb{C}\setminus(-\infty,0]$. A residue calculation yields that

$$h(z)=\frac{k^2}{2zR(z)}\log\left(\frac{\sqrt{ab}+z-R(z)}{\sqrt{ab}+z+R(z)}\right),\qquad z\in\mathbb{C}\setminus(-\infty,0], \tag{92}$$

where log denotes the principal branch of logarithm. For a computation below, we note that as $k\to\infty$,

$$\sqrt{\mathbf{a}}=1-\sqrt{\frac{2}{k}}+\frac{2}{k}+O(k^{-3/2}),$$
$$\sqrt{\mathbf{b}}=1+\sqrt{\frac{2}{k}}+\frac{2}{k}+O(k^{-3/2}). \tag{93}$$

We also remark that with $x = 1 + \frac{2w}{\sqrt{k}}$, for $w = O(1)$, as $k \to \infty$, at least formally,

$$\psi(x)\,dx \sim \frac{k}{\pi}\sqrt{2 - w^2}\,dw, \qquad w \in [-\sqrt{2}, \sqrt{2}], \tag{94}$$

which is precisely Wigner's semicircle. This last calculation is not going to be used below, but it provides an intuitive reason as to why the ensemble (74) [and (72)] has the same asymptotics as the Gaussian unitary ensemble, not only locally, but also globally.

Set

$$\mathbf{M}(z) = e^{-(1/2)\ell\sigma_3}\mathbf{Y}(z)e^{-G(z)\sigma_3}e^{(1/2)\ell\sigma_3} \tag{95}$$

for $z \in \mathbb{C} \setminus \overline{\mathbb{R}}_+$. Using the analyticity of $G$ for $z \in \mathbb{R}_+ \setminus [\mathbf{a}, \mathbf{b}]$ and the variational conditions, $\mathbf{M}(z)$ solves the following, equivalent, RHP:

- $\mathbf{M}(z)$ is analytic for $z \in \mathbb{C} \setminus \overline{\mathbb{R}}_+$, $\mathbf{M}_\pm(z)$ is continuous for $z \in \mathbb{R}_+$ and $\mathbf{M}(z)$ is bounded as $z \to 0$;
- for $z \in \mathbb{R}_+$, $\mathbf{M}_+(z) = \mathbf{M}_-(z)\mathbf{V}_M(z)$, where

$$\mathbf{V}_M(z) = \begin{pmatrix} e^{G_-(z)-G_+(z)} & 1 \\ 0 & e^{G_+(z)-G_-(z)} \end{pmatrix}, \qquad z \in (\mathbf{a}, \mathbf{b}), \tag{96}$$

$$\mathbf{V}_M(z) = \begin{pmatrix} 1 & e^{2G(z)+\log(w(z))-\ell} \\ 0 & 1 \end{pmatrix}, \qquad z \in \mathbb{R}_+ \setminus (\mathbf{a}, \mathbf{b}); \tag{97}$$

- $\mathbf{M}(z) = \mathbf{I} + O(z^{-1})$ as $z \to \infty$.

The nonunit terms in the jump matrix can be expressed in a unifying way. Set

$$H(z) = G(z) + \tfrac{1}{2}\log(w(z)) - \tfrac{1}{2}\ell, \qquad z \in \mathbb{C} \setminus ((-\infty, 0] \cup [\mathbf{a}, \mathbf{b}]). \tag{98}$$

Noting the variational condition, we find that for $z \in (\mathbf{a}, \mathbf{b})$,

$$\begin{aligned} G_+(z) - G_-(z) &= 2G_+(z) + \log(w(z)) - \ell = 2H_+(z) \\ &= -(2G_-(z) + \log(w(z)) - \ell) = -2H_-(z). \end{aligned} \tag{99}$$

Hence, $G_+(z) - G_-(z)$ has an analytic continuation both above and the below the real axis. Therefore, the jump matrix $\mathbf{V}_M$ equals

$$\mathbf{V}_M(z) = \begin{pmatrix} e^{-2H_+(z)} & 1 \\ 0 & e^{-2H_-(z)} \end{pmatrix}, \qquad z \in (\mathbf{a}, \mathbf{b}), \tag{100}$$

$$\mathbf{V}_M(z) = \begin{pmatrix} 1 & e^{2H(z)} \\ 0 & 1 \end{pmatrix}, \qquad z \in \mathbb{R}_+ \setminus (\mathbf{a}, \mathbf{b}). \tag{101}$$

Using the definition of $G$ and Lemma 4, one can check that

$$H'(z) = \tfrac{1}{2}R(z)h(z). \tag{102}$$

We now scale the RHP for $\mathbf{M}$ so that the interval $(\mathbf{a},\mathbf{b})$ becomes $(-1,1)$. In other words, instead of moving the interval as the support of the equilibrium measure, we will fix the support. In that way, we can use the analysis of [20, 21] more directly. Define

$$\mathbf{N}(z)=\mathbf{M}\Big(\frac{\mathbf{b}-\mathbf{a}}{2}z+\frac{\mathbf{b}+\mathbf{a}}{2}\Big). \tag{103}$$

Set $\Sigma=(-\frac{\mathbf{b}+\mathbf{a}}{\mathbf{b}-\mathbf{a}},\infty)$ and set

$$\hat{H}(z)=H\Big(\frac{\mathbf{b}-\mathbf{a}}{2}z+\frac{\mathbf{b}+\mathbf{a}}{2}\Big). \tag{104}$$

The matrix $\mathbf{N}$ solves the following RHP:

- $\mathbf{N}(z)$ is analytic for $z\in\mathbb{C}\setminus\overline{\Sigma}$, $\mathbf{N}_\pm(z)$ is continuous for $z\in\Sigma$ and $\mathbf{N}(z)$ is bounded as $z\to-\frac{\mathbf{b}+\mathbf{a}}{\mathbf{b}-\mathbf{a}}$;
- for $z\in\Sigma$, $\mathbf{N}_+(z)=\mathbf{N}_-(z)\mathbf{V}_N(z)$, where

$$\mathbf{V}_N(z)=\begin{pmatrix} e^{-2\hat{H}_+(z)} & 1\\ 0 & e^{-2\hat{H}_-(z)}\end{pmatrix}, \qquad z\in(-1,1), \tag{105}$$

$$\mathbf{V}_N(z)=\begin{pmatrix} 1 & e^{2\hat{H}(z)}\\ 0 & 1\end{pmatrix}, \qquad z\in\Sigma\setminus(-1,1); \tag{106}$$

- $\mathbf{N}(z)=\mathbf{I}+O(z^{-1})$ as $z\to\infty$.

Note the factorization for $z\in(-1,1)$,

$$\begin{pmatrix} e^{-2\hat{H}_+(z)} & 1\\ 0 & e^{-2\hat{H}_-(z)}\end{pmatrix}=\begin{pmatrix} 1 & 0\\ e^{-2\hat{H}_-(z)} & 1\end{pmatrix}\begin{pmatrix} 0 & 1\\ -1 & 0\end{pmatrix}\begin{pmatrix} 1 & 0\\ e^{-2\hat{H}_+(z)} & 1\end{pmatrix}, \tag{107}$$

where we use the fact that $\hat{H}_+(z)+\hat{H}_-(z)=0$ for $z\in(-1,1)$. Let $\Sigma_j$, $j=0,1,\dots,4$, and $\Omega_j$, $j=1,\dots,4$, be the contours and open regions given in Figure 1. Contours are oriented from left to right. Define

$$\mathbf{Q}(z)=\begin{cases}\mathbf{N}(z), & z\in\Omega_1\cup\Omega_4,\\ \mathbf{N}(z)\begin{pmatrix} 1 & 0\\ -e^{-2\hat{H}(z)} & 1\end{pmatrix}, & z\in\Omega_2,\\ \mathbf{N}(z)\begin{pmatrix} 1 & 0\\ e^{-2\hat{H}(z)} & 1\end{pmatrix}, & z\in\Omega_3.\end{cases} \tag{108}$$

Then $\mathbf{Q}_+(z)=\mathbf{Q}_-(z)\mathbf{V}_Q(z)$ for $z$ in $\Sigma_0,\dots,\Sigma_4$, where

$$\mathbf{V}_Q(z)=\begin{pmatrix} 0 & 1\\ -1 & 0\end{pmatrix}, \qquad z\in\Sigma_0, \tag{109}$$

$$\mathbf{V}_Q(z)=\begin{pmatrix} 1 & 0\\ e^{-2\hat{H}(z)} & 1\end{pmatrix}, \qquad z\in\Sigma_1\cup\Sigma_2, \tag{110}$$

$$\mathbf{V}_Q(z)=\begin{pmatrix} 1 & e^{2\hat{H}(z)}\\ 0 & 1\end{pmatrix}, \qquad z\in\Sigma_3\cup\Sigma_4. \tag{111}$$

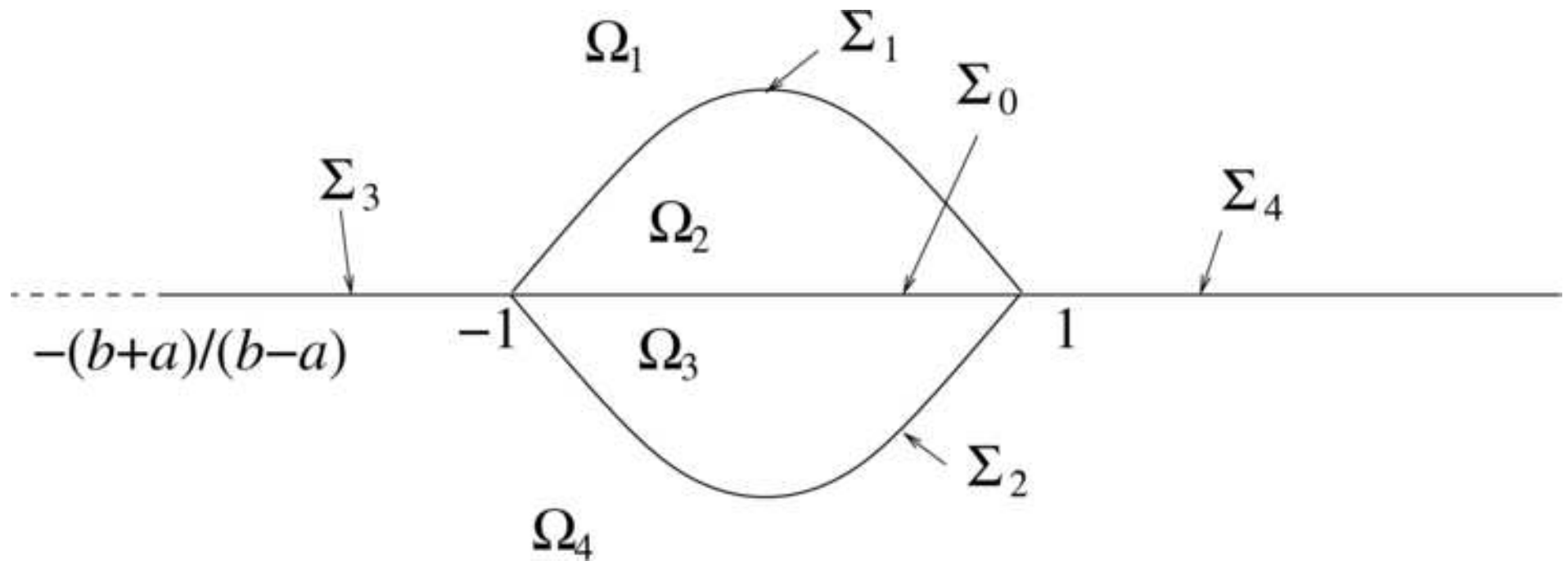

FIG. 1. *Contours for* $\mathbf{N}$.

The off-diagonal terms of $\mathbf{V}_Q$ on $\Sigma_1 \cup \cdots \cup \Sigma_4$ converge to 0 as the following lemma implies.

LEMMA 5. *There exist* $\delta_0 > 0$ *and* $k_0 > 0$ *such that for* $k \geq k_0$,

$$\text{Re}[\hat{H}(x+iy)] \geq 2k|y|\sqrt{1-x^2} \qquad \text{for } -1 \leq x \leq 1 \text{ and } -\delta_0 \leq y \leq \delta_0. \tag{112}$$

*For any* $\delta > 0$,

$$\hat{H}(x) \leq -k\delta^{3/2} \qquad \text{for } -\frac{\mathbf{b}+\mathbf{a}}{\mathbf{b}-\mathbf{a}} < x \leq -1-\delta \text{ and } x \geq 1+\delta \tag{113}$$

*when* $k \geq k_0$, *and*

$$\lim_{k\to\infty} \int_{(\Sigma_3\cup\Sigma_4)\cap\{|z-1|>\delta\}\cap\{|z+1|>\delta\}} e^{2\hat{H}(x)}\,dx = 0. \tag{114}$$

Hence, $\mathbf{V}_Q \to \mathbf{V}_\infty$ for a constant matrix $\mathbf{V}_\infty$ defined as

$$\mathbf{V}_\infty(z) = \begin{pmatrix} 0 & 1 \\ -1 & 0 \end{pmatrix}, \qquad z \in \Sigma_0, \tag{115}$$

and $\mathbf{V}_\infty(z) = \mathbf{I}$ for $z \in \Sigma_1 \cup \cdots \cup \Sigma_4$, where the convergence $\mathbf{V}_Q \to \mathbf{V}_\infty$ is in $L^\infty(\Sigma_0 \cup \cdots \cup \Sigma_4)$ and also in $L^2((\Sigma_0 \cup \cdots \cup \Sigma_4) \cap \{|z-1| > \delta\} \cap \{|z+1| > \delta\})$ for an arbitrary, but fixed, $\delta > 0$. Let

$$\beta(z) = \left(\frac{z-1}{z+1}\right)^{1/4}, \tag{116}$$

where the branch cut is $[-1,1]$ and $\beta(z) \sim 1$ as $z \to +\infty$ on the real line, and define

$$\mathbf{Q}^\infty(z) = \frac{1}{2}\begin{pmatrix} \beta+\beta^{-1} & -i(\beta-\beta^{-1}) \\ i(\beta-\beta^{-1}) & \beta+\beta^{-1} \end{pmatrix} \tag{117}$$

for $z \in \mathbb{C} \setminus \Sigma_0$. Then $\mathbf{Q}^\infty(z)$ is the solution to the RHP for the $\mathbf{Q}^\infty_+ = \mathbf{Q}^\infty_- \mathbf{V}_\infty$ and $\mathbf{Q}^\infty(z) \to \mathbf{I}$ as $z \to \infty$. The convergence $\mathbf{V}_Q \to \mathbf{V}_\infty$ is not uniform near the points $z = \pm 1$, hence it is not true that $\mathbf{Q}(z) \to \mathbf{Q}^{\infty}(z)$ for all $z$ and one therefore needs local parametrix for $z$ in a neighborhood of $\pm 1$.

Let $\Psi(z)$ be the matrix-valued function constructed from the Airy function and its derivatives, as defined in Proposition 7.3 of [20]. Let $\varepsilon > 0$. For $z \in U_r := \{z : |z-1| < \varepsilon\}$, set

$$\mathbf{S}_r(z) = \mathbf{E}(z)\Psi((-\tfrac{3}{2}\hat{H}(z))^{2/3})e^{-\hat{H}(z)\sigma_3}, \tag{118}$$

where

$$\begin{aligned}\mathbf{E}(z) = {}& \sqrt{\pi}e^{(\pi/6)i}\begin{pmatrix} 1 & -1 \\ -i & -i \end{pmatrix} \\ & \times \begin{pmatrix} (-\frac{3}{2}\hat{H}(z))^{1/6}\beta(z)^{-1} & 0 \\ 0 & (-\frac{3}{2}\hat{H}(z))^{-1/6}\beta(z) \end{pmatrix}.\end{aligned} \tag{119}$$

Note that $\mathbf{E}(z)$ is analytic in $U_r$ if $\varepsilon$ is chosen sufficiently small. The matrix $\mathbf{S}_l(z)$ is defined in a similar way for $z \in U_l := \{z : |z+1| < \varepsilon\}$. Define

$$\mathbf{Q}_{\text{par}}(z) = \begin{cases} \mathbf{Q}^{\infty}(z), & z \in \mathbb{C} \setminus U_r \cup U_l \cup \Sigma, \\ \mathbf{S}_r(z), & z \in U_r \setminus \Sigma, \\ \mathbf{S}_l(z), & z \in U_l \setminus \Sigma. \end{cases} \tag{120}$$

From the basic theory of RHP, the estimate in Lemma 5 and the same argument as in [20], one can check that the jump matrix for $\mathbf{Q}_{\text{par}}^{-1}\mathbf{Q}$ converges to the identity in $L^2 \cap L^{\infty}$. Hence,

$$\mathbf{Q}(z) = (\mathbf{I} + O(k^{-1}))\mathbf{Q}_{\text{par}}(z). \tag{121}$$

This holds uniformly for $z$ outside an open neighborhood of the contours $\Sigma \cup \partial U_r \cup \partial U_l$. But a simple deformation argument implies that the result is extended to $z$ on the contours (see [20]). Hence, by reversing the transformations $\mathbf{Y} \to \mathbf{M} \to \mathbf{N} \to \mathbf{Q}$ [see (95), (103) and (108)], the asymptotics of $\mathbf{Y}(z)$ for all $z \in \mathbb{C}$ are obtained.

By substituting the asymptotics of $\mathbf{Y}$ into (84), edge and bulk scaling limits of the $\mathbf{K}_k$ are obtained; see [15, 18, 21] for details. For $x_0$ such that $\sqrt{k}(x_0 - 1)$ lies in a compact subset of $(\sqrt{k}(\mathbf{a}-1), \sqrt{k}(\mathbf{b}-1))$, for all $\xi, \eta$ in a compact subset of $\mathbb{R}$,

$$\frac{1}{\psi(x_0)}\mathbf{K}_k\left(x_0 + \frac{\xi}{\psi(x)}, x + \frac{\eta}{\psi(x_0)}\right) \to \mathbb{S}(\xi, \eta) \tag{122}$$

in trace norm for $\xi, \eta \in \mathbb{R}$, where

$$\mathbb{S}(\xi, \eta) = \frac{\sin(\pi(\xi - \eta))}{\pi(\xi - \eta)}. \tag{123}$$

Here, we may replace $\psi(x_0)$ by $\mathbf{K}_k(x_0, x_0)$. The error is $O(k^{-1})$, uniformly for $\xi, \eta$ in a compact set. The convergence is also in trace norm in the Hilbert space $L^2((-\eta,\eta))$ for a fixed $\eta > 0$. From (82), by taking $x_0 = e^{2/k}$, the limit (30) in Proposition 2 is obtained.

At the edge of the support of $\psi(x)$, set

$$B_k = \left[-\frac{1}{2}\sqrt{\mathbf{b}-\mathbf{a}}h(b)\right]^{2/3} \sim \frac{k^{7/6}}{\sqrt{2}}. \tag{124}$$

As $k \to \infty$,

$$\frac{1}{B_k}\mathbf{K}_k\left(b+\frac{\xi}{B_k}, b+\frac{\eta}{B_k}\right) \to \mathbb{A}(\xi,\eta) \tag{125}$$

in trace norm in the Hilbert space $L^2((\xi,\infty))$ for a fixed $\xi$, where

$$\mathbb{A}(\xi,\eta) = \frac{\mathrm{Ai}(\xi)\mathrm{Ai}'(\eta) - \mathrm{Ai}'(\xi)\mathrm{Ai}(\eta)}{\xi-\eta} \tag{126}$$

is the Airy kernel. Hence, from (82), the limit (31) in Proposition 2 is obtained.

**4. Generalizations and discussions.** We comment on three issues in this section: the case in which the moment generating function does not exist, finite-dimensional distributions and the connection of this work to $q$-orthogonal polynomials.

*No moment generating function.* In this paper, we have assumed the existence of the moment generating function for the random variable increments of nonintersecting random walks. This is simply to improve the estimates. For the case $\mathbb{E}|X_i^j|^{2+\delta} < \infty$, $\delta > 0$, there is a version of the KMT theorem which provides analogous estimates to those used in Section 2. As one would expect for this case, $N_k$ must grow more rapidly in $k$. Another method for achieving results similar to those of this paper is to use Skorohod embedding in order to embed the nonintersecting random walks into Brownian motions. In order to achieve this, one must assume that $\mathbb{E}|X_i^j|^4 < \infty$.

*Finite-dimensional distributions.* The results of this paper focus on the limiting distributions of nonintersecting random walks at the fixed time $t = 1$. It is also interesting to consider finite-dimensional distributions of the process, that is, in the correct scaling $t_1, \dots, t_n \in [1 - Ak^{-1/3}, 1 + Ak^{-1/3}]$, the finite-dimensional distributions of the fluctuations of the top random walk should converge to those of the Airy process. A similar, but differently scaled, result should also be true "in bulk"; see, for example, [1, 46, 53] and references therein concerning the Airy process and other processes from random matrix theory. The methods of Section 2 are certainly applicable to this problem, however, the convergence of the finite-dimensional distributions of the nonintersecting Brownian bridges to Airy/sine processes does not follow immediately from the analysis of Section 3. However, one can use a different approach based on the method of Eynard and Mehta [14, 23]. In this approach, an inversion of a matrix is crucial. After the completion of the present paper, Widom communicated to the authors how to invert the matrix. Work in this direction will appear in a future paper.

*Stieltjes–Wigert weight and $q$-orthogonal polynomials.* In Section 3, the Riemann–Hilbert problems for the orthogonal polynomials with respect to the Stieltjes–Wigert weight (79) was analyzed in the Plancherel–Rotach asymptotic regime. The analysis yields the asymptotics of the Stieltjes–Wigert polynomials in the entire complex plane. Since Stieltjes–Wigert polynomials are examples of $q$-polynomials, this result also yields an asymptotic result for certain $q$-polynomials.

**Acknowledgments.** The authors would like to thank Percy Deift and Harold Widom for useful discussions.

Department of Mathematics
University of Michigan
Ann Arbor, Michigan 48109
USA
E-mail: baik@umich.edu

Department of Mathematics
University of California
Santa Cruz, California 95064
USA
E-mail: tsuidan@ucsc.edu